\documentclass[10pt,sumlimits,namelimits]{article}
\newcommand{\nc}{\newcommand}
\newcommand{\rnc}{\renewcommand}
\usepackage{amsthm}
\usepackage{amsmath}
\usepackage{amsfonts}
\usepackage{amssymb}
\usepackage{color}
\usepackage[colorlinks=true,urlcolor=blue]{hyperref}
\usepackage[applemac]{inputenc}

\nc{\N}{\mathbb{N}}
\nc{\Z}{\mathbb{Z}}
\nc{\D}{\mathbb{D}}
\nc{\Q}{\mathbb{Q}}
\nc{\R}{\mathbb{R}}
\nc{\C}{\mathbb{C}}
\nc{\vD}{{\cal D}}
\nc{\vS}{{\cal S}}
\nc{\vphi}{\varphi}
\nc{\eps}{\varepsilon}
\nc{\dsp}{\displaystyle}
\nc{\ovl}{\overline}
\nc{\udl}{\underline}
\nc{\vlim}{\lim\limits}
\nc{\vlimsup}{\limsup\limits}
\nc{\vliminf}{\liminf\limits}
\nc{\vsup}{\sup\limits}
\nc{\vinf}{\inf\limits}
\nc{\vint}{\int\limits}
\nc{\inj}{\hookrightarrow}
\nc{\tends}{\longrightarrow}
\nc{\weak}{\rightharpoonup}
\nc{\w}{{\textsl w}}
\nc{\loc}{{\rm loc}}
\rnc{\le}{\leqslant}
\rnc{\ge}{\geqslant}
\rnc{\Re}{{\rm Re}}
\rnc{\Im}{{\rm Im}}

\numberwithin{equation}{section}
\rnc{\theequation}{\thesection.\arabic{equation}}

\newtheorem{thm}{Theorem}[section]
\newtheorem{prop}[thm]{Proposition}
\newtheorem{cor}[thm]{Corollary}
\newtheorem{lem}[thm]{Lemma}

\theoremstyle{definition}
\newtheorem{rmk}[thm]{Remark}
\newtheorem{defi}[thm]{Definition}

\newenvironment{proof*}{\noindent{\bf Proof.}}{\qed}
\newenvironment{vproof}[1]{\noindent{\bf Proof #1}}{\qed}

\title{\Huge \sc Maximum Decay Rate for the Nonlinear Schrödinger Equation}
\author{\sc Pascal Bégout}
\date{}

\begin{document}

\maketitle

\begin{center}
Laboratoire Jacques-Louis Lions \\
Université Pierre et Marie Curie \\
Boîte Courrier 187 \\
4, place Jussieu 75252 Paris Cedex 05, FRANCE  \bigskip \\
{\footnotesize e-mail\:: }\htmladdnormallink{{\footnotesize\udl{\tt{begout@ann.jussieu.fr}}}}
{mailto:begout@ann.jussieu.fr}
\end{center}

\begin{abstract}
In this paper, we consider global solutions for the following nonlinear Schrödinger equation $iu_t+\Delta u+\lambda|u|^\alpha u=0,$ in $\R^N,$ with
$\lambda\in\R$ and $0\le\alpha<\frac{4}{N-2}$ $(0\le\alpha<\infty$ if $N=1).$ We show that no nontrivial solution can decay faster than the solutions of the free Schrödinger equation, provided that $u(0)$ lies in the weighted Sobolev space $H^1(\R^N)\cap L^2(|x|^2;dx),$ in the energy space, namely $H^1(\R^N),$ or in $L^2(\R^N),$ according to the different cases.
\end{abstract}

\baselineskip .7cm

\section{Introduction and notations}
\label{introduction}

{\let\thefootnote\relax\footnotetext{2000 Mathematics Subject Classification: 35Q55,(35B40)}}

We consider global solutions of the following nonlinear Schrödinger equation,
\begin{gather}
 \left\{
  \begin{split}
   \label{nls}
    i\frac{\partial u}{\partial t}+\Delta u+\lambda|u|^\alpha u & = 0,\; (t,x)\in[0,\infty)\times \R^N, \\
                                                           u(0) & = \vphi, \mbox{ in }\R^N,
  \end{split}
 \right.
\end{gather}
where $\lambda\in \R,$ $0\le\alpha<\dfrac{4}{N-2}$ $(0\le\alpha<\infty$ if $N=1)$ and $\vphi$ a given initial data.

It is well-known that if we denote by $T(t)$ the Schrödinger's free operator, then for every $r\in[2,\infty]$ and for
every $\vphi\in L^{r'}(\R^N),$
\begin{gather}
 \label{groupe}
  \forall t\in\R\setminus\{0\},\; \|T(t)\vphi\|_{L^r}\le (4\pi|t|)^{-N\left(\frac{1}{2}-\frac{1}{r}\right)}\|\vphi\|_{L^{r'}},
\end{gather}
where $r'=\frac{r}{r-1}.$ Note that for every $r\in\left[2,\frac{2N}{N-2}\right)$ $(r\in[2,\infty)$ if $N=1),$ if $\vphi\in
L^2(\R^N)\cap L^2(|x|^2;dx)$ then $\vphi\in L^{r'}(\R^N)$ and $\|\vphi\|_{L^{r'}}\le C(\|\vphi\|_{L^2},\|x\vphi\|_{L^2}).$
Furthermore, the estimate (\ref{groupe}) is optimal in the following sense. For every $r\in[1,\infty]$ and for every $\vphi\in\vS'(\R^N),$ if $\vphi\not\equiv0$ then
\begin{gather}
 \label{maxdecaylinear}
  \liminf_{t\to\pm\infty}|t|^{N\left(\frac{1}{2}-\frac{1}{r}\right)}\|T(t)\vphi\|_{L^r}>0.
\end{gather}
For the proof, see Strauss \cite{297.35062} $(r\ge2)$ and Kato \cite{MR95i:35276} (general case). In the same way, there exist some solutions of the
nonlinear Schrödinger equation (\ref{nls}) which have a linear decay (in the sense of (\ref{groupe})). See for example Cazenave
\cite{caz1}, Theorem 7.2.1; Hayashi and Naumkin \cite{MR99f:35190}. In particular, these solutions lie in $H^1(\R^N)\cap
L^2(|x|^2;dx).$ On the other hand, we know that there exist solutions of some heat, Ginzburg-Landau and Schrödinger type equations
which have a decay rate faster than the corresponding linear problem (Hayashi, Kaikina and Naumkin \cite{MR1785915,MR1800091}).
Take an example. Let $u$ be a classical solution of the heat equation $u_t-\Delta u+|u|^\alpha u=0,$
$(t,x)\in[0,\infty)\times\R^N,$ with initial datum $u_0\in L^\infty(\R^N),$ $u_0\not\equiv0$ and $u_0\ge0$ a.e. Then we have
by the maximum principle, $\|u(t)\|_{L^\infty}\le(\alpha t)^{-\frac{1}{\alpha}},$ for every $t>0,$ whereas for every $t\ge1,$
$\|e^{t\Delta}u_0\|_{L^\infty}\ge Ct^{-\frac{N}{2}},$ for some constant $C>0.$ Thus, if $0<\alpha<\frac{2}{N}$ then $u(t)$ decays
faster than $e^{t\Delta}u_0.$ So, it is natural to wonder if some solutions of the nonlinear Schrödinger equation (\ref{nls}) may
have faster decay than the solutions of the linear equation. We will see that such solutions do not exist (except the trivial
solution). There exist partial results in this direction. This is the case for $\alpha=\frac{2}{N}$ and $N=1$ (Hayashi and Naumkin
\cite{MR99f:35190}), for $\alpha=\frac{4}{N}$ (Cazenave and Weissler \cite{MR92c:35111}, Theorem 2.1 (a)) or  for some
self-similar solutions for $\alpha>\alpha_0,$ where $\alpha_0$ is given by $\alpha_0=\frac{-(N-2)+\sqrt{N^2+12N+4}}{2N}$
(Cazenave and Weissler \cite{MR1745480}, Corollary 3.9).

This paper is organized as follows. In Section \ref{decmaxresprincX}, we give the main results concerning the solutions lying in
$H^1(\R^N)\cap L^2(|x|^2;dx)$ and in $H^1(\R^N).$ In Section \ref{decmaxresprincL2}, we give the main results concerning the
solutions lying in $L^2(\R^N).$ In Section \ref{decmaxestiminfini}, we give several estimates for large times and establish Lemma
\ref{lemmaxdecaynonlinear3}, which asserts that the existence of a scattering state in $L^2(\R^N)$ implies a maximum rate decay
which is linear (in the sense that the solution satisfies (\ref{maxdecaylinear})). Lemma \ref{lemmaxdecaynonlinear3} is at the
heart of the results of this paper. Finally, we will prove the results for solutions in $H^1(\R^N)\cap L^2(|x|^2;dx)$ and
$H^1(\R^N)$ in Section \ref{vproofdecmaxresprincX}, and those for solutions in $L^2(\R^N)$ in Section \ref{vproofdecmaxresprincL2}.

We finish this section by giving some notations and we recall an embedding property of the weighted Sobolev space $L^2(\R^N)\cap
L^2(|x|^2;dx),$ which will be used to prove the results for solutions lying in this space, and some results of the solutions of
the nonlinear Schrödinger equation (\ref{nls}).

We design by $\ovl{z}$ the conjugate of the complex number $z$ and
$\Delta=\sum\limits_{j=1}^N\frac{\partial^2}{\partial x_j^2}.$
For $p\in[1,\infty],$ we denote by $p'$ the conjugate of $p$ defined by $\frac{1}{p}+\frac{1}{p'}=1$ and by
$L^p(\R^N) = L^p(\R^N;\C),$ with norm $\|\: .\:\|_{L^p},$ the Lebesgue spaces.
$H^1(\R^N) = H^1(\R^N;\C)$ with norm $\|\: .\:\|_{H^1},$ is the well-known Sobolev space and we use the convention
$W^{0,p}(\R^N)=L^p(\R^N)$ and $H^0(\R^N)=W^{0,2}(\R^N)=L^2(\R^N).$
We define the Hilbert spaces $Y=\left\{\psi\in L^2(\R^N;\C);\; \|\psi\|_Y<\infty\right\}$ with norm
$\|\psi\|_Y^2=\|\psi\|_{L^2(\R^N)}^2+\vint_{\R^N}|x|^2|\psi(x)|^2dx$ and
$X=\left\{\psi\in H^1(\R^N;\C);\; \|\psi\|_X<\infty\right\}$ with norm $\|\psi\|_X^2=\|\psi\|_{H^1(\R^N)}^2+
\vint_{\R^N}|x|^2|\psi(x)|^2dx.$
For a functional space $E\subset\vS'(\R^N)$ with norm $\| \: . \: \|_E,$ we write $\|f\|_E=\infty$ if $f\in\vS'(\R^N)$ and
if $f\not\in E.$ We design by $(T(t))_{t \in \R}$ the group of isometries $(e^{it\Delta})_{t\in\R}$ generated by $i\Delta$ on
$L^2(\R^N;\C)$ and by $C$ the auxiliary positive constants.
Finally $C(a_1, a_2,\dots,a_n)$ indicates that the constant $C$ depends only on parameters $\:a_1,a_2,\dots, a_n$ and that the
dependence is continuous.

It is clear that $Y$ is a separable Hilbert space and that $Y\inj L^{r'}(\R^N)$ with dense embedding, if
$r\in\left[2,\frac{2N}{N-2}\right)$ $(r\in[2,\infty]$ if $N=1).$

We recall that for every $\vphi\in H^1(\R^N),$ (\ref{nls}) has a unique solution $u\in C((-T_*,T^*);H^1(\R^N))$
which satisfies the conservation of charge and energy, that is, for all $t\in(-T_*,T^*),$ $\|u(t)\|_{L^2}=\|\vphi\|_{L^2}$
and $E(u(t))=E(\vphi),$ where $E(\vphi)\stackrel{\rm def}{=}\frac{1}{2}\|\nabla\vphi\|_{L^2}^2
-\frac{\lambda}{\alpha+2}\|\vphi\|_{L^{\alpha+2}}^{\alpha+2}.$ Moreover, for every admissible pair $(q,r)$ (see Definition
\ref{padef} below), $u\in L^q_\loc((-T_*,T^*);W^{1,r}(\R^N))).$ In addition, if $\lambda\le0,$ if $\alpha<\frac{4}{N}$ or if
$\|\vphi\|_{H^1}$ is small enough then $T^*=T_*=\infty$ and $\|u\|_{L^\infty(\R;H^1)}<\infty.$ Finally, if $\vphi\in X$ then
$u\in C((-T_*,T^*);X).$ See Ginibre and Velo \cite{MR82c:35059,MR82c:35057,MR82c:35058,MR87a:35164,MR87b:35150}, Kato
\cite{MR88f:35133,MR91d:35202}. See also Cazenave and Weissler \cite{MR89j:35114,MR91c:35136}. We are also interested by solutions
in $L^2(\R^N).$ We recall that if $0<\alpha\le\frac{4}{N}$ then for every $\vphi\in L^2(\R^N),$ (\ref{nls}) has a unique solution
$u\in C((-T_*,T^*);L^2(\R^N))\cap L^q_\loc((-T_*,T^*);L^{\alpha+2}(\R^N)),$ where $q=\frac{4(\alpha+2)}{N\alpha},$ which satisfies
the above conservation of charge. In addition, for every admissible pair $(q,r),$ $u\in L^q_\loc((-T_*,T^*);L^r(\R^N)).$ Finally,
if $\alpha<\frac{4}{N}$ then $T^*=T_*=\infty.$ See Tsutsumi \cite{MR89c:35143}. See also Cazenave and Weissler
\cite{MR91a:35149,MR91j:35252}.

\begin{defi}
\label{padef}
We say that $(q,r)$ is an {\it admissible pair} if the following holds. \medskip \\
$
\begin{array}{rl}
 (\rm i) & 2\le r\le\frac{2N}{N-2}\; (2\le r<\infty\; $if$\; N=2,\; 2\le r\le\infty\mbox\; $if$\; N=1), \medskip \\
(\rm ii) & \frac{2}{q} = N\left(\frac{1}{2} - \frac{1}{r}\right).
\end{array}
$
\\
Note that in this case $2\le q\le\infty$ and $q=\dfrac{4r}{N(r-2)}.$
\end{defi}

Finally, we recall the Strichartz' estimates. Let $I\subseteq\R,$ be an interval, let $t_0\in\ovl I,$ let $(q,r)$ and
$(\gamma,\rho)$ be two admissible pairs, let $\vphi\in L^2(\R^N)$ and let $f\in L^{\gamma'}(I;L^{\rho'}(\R^N)).$ Then the
following integral equation defined for all $t\in I,$ $u(t)=T(t)\vphi+i\dsp\int_{t_0}^tT(t-s)f(s)ds,$ satisfies the following
inequality
$
\|u\|_{L^q(I,L^r)}\le C_0\|\vphi\|_{L^2}+C_1\|f\|_{L^{\gamma'}(I;L^{\rho'})},
$
where $C_0=C_0(N,r)$ and $C_1=C_1(N,r,\rho).$ For more details, see Keel and Tao \cite{MR1646048}.

\section{Sharp lower bound}
\label{decmaxresprincX}

\begin{thm}
\label{maxdecaynonlinearX-I}
Let $\lambda\le0,$ $0\le\alpha<\dfrac{4}{N-2}$ $(0\le\alpha<\infty$ if $N=1),$ $\vphi\in H^1(\R^N)$ and $u$ be the corresponding
solution of $(\ref{nls}).$ If $\alpha\le\dfrac{4}{N}$ then assume further that $\vphi\in X.$ If $\vphi\not\equiv0$ then for every
$r\in[2,\infty],$
\begin{gather*}
\vliminf_{t\to\infty}|t|^{N\left(\frac{1}{2}-\frac{1}{r}\right)}\|u(t)\|_{L^r}>0.
\end{gather*}
\end{thm}

\begin{thm}
\label{maxdecaynonlinearX-II}
Let $\lambda>0,$ $0\le\alpha<\dfrac{4}{N-2}$ $(0\le\alpha<\infty$ if $N=1)$ and $\vphi\in H^1(\R^N)$ be such that the
corresponding solution $u$ of $(\ref{nls})$ is positively global in time. If $\alpha\le\dfrac{4}{N}$ then assume further that
$\vphi\in X.$ If $\vphi\not\equiv0$ then for every $r\in[\alpha+2,\infty],$
$$
\begin{array}{rl}
 & \left\{
  \begin{array}{rl}
   \vliminf_{t\to\infty}|t|^{N\left(\frac{1}{2}-\frac{1}{r}\right)}\|u(t)\|_{L^r}>0, & if\;\; \alpha\le\dfrac{4}{N}, \medskip \\
   \vlimsup_{t\to\infty}|t|^{N\left(\frac{1}{2}-\frac{1}{r}\right)}\|u(t)\|_{L^r}>0, & if \;\;\alpha > \dfrac{4}{N}.
  \end{array}
 \right.
\end{array}
$$
And if there exists $\rho\in[\alpha+2,\infty]$ such that
$\vlimsup_{t\to\infty}|t|^{N\left(\frac{1}{2}-\frac{1}{\rho}\right)}\|u(t)\|_{L^\rho}<\infty$ then for every $r\in[2,\infty],$
$$
\vliminf_{t\to\infty}|t|^{N\left(\frac{1}{2}-\frac{1}{r}\right)}\|u(t)\|_{L^r}>0.
$$
\end{thm}

\begin{rmk}
Theorems \ref{maxdecaynonlinearX-I} and \ref{maxdecaynonlinearX-II} assert that if $u$ is a solution of (\ref{nls}) with
$\lambda\in\R,$ $\alpha=0$ and initial data $\vphi\in X,$ then for every $r\in[2,\infty],$
$\vliminf_{t\to\infty}|t|^{N\left(\frac{1}{2}-\frac{1}{r}\right)}\|u(t)\|_{L^r}>0.$
\end{rmk}

In the attractive case and when $\alpha>\frac{4}{N+2}$ $(\alpha>2$ if $N=1),$ we may obtain an optimal lower bound.
It is sufficient to choose $\|\vphi\|_X$ small enough (see corollary below).

\begin{cor}
\label{cormaxdecaynonlinearX-II}
Let $\lambda>0,$ $\dfrac{4}{N+2}<\alpha<\dfrac{4}{N-2}$ $(2<\alpha<\infty$ if $N=1),$ $\vphi\in X$ and $u$ be the
corresponding solution $u$ of $(\ref{nls}).$ If $\vphi\not\equiv0$ and if $\|\vphi\|_X$ is small enough then
$u$ is global in time and for every $r\in[2,\infty],$
\begin{gather*}
\vliminf_{t\to\infty}|t|^{N\left(\frac{1}{2}-\frac{1}{r}\right)}\|u(t)\|_{L^r}>0.
\end{gather*}
\end{cor}

When $\alpha=\frac{4}{N},$ we may suppose that $\vphi\in H^1(\R^N)$ instead of $\vphi\in X,$ as shows the following proposition,
provided that $\|\vphi\|_{H^1}$ is small enough.

\begin{prop}
\label{maxdecaynonlinearX-III}
Let $\lambda\in\R\setminus\{0\},$ $\alpha=\dfrac{4}{N},$ $\vphi\in H^1(\R^N)$ and $u$ be the associated solution of $(\ref{nls}).$
If $\vphi\not\equiv0$ and if $\|\vphi\|_{H^1}$ is small enough then $u$ is global in time and for every $r\in[2,\infty],$
\begin{gather*}
\vliminf_{t\to\infty}|t|^{N\left(\frac{1}{2}-\frac{1}{r}\right)}\|u(t)\|_{L^r}>0.
\end{gather*}
\end{prop}

\begin{rmk}
\label{rmkmaxdecaynonlinearX-III}
The lower bounds obtained in Theorems \ref{maxdecaynonlinearX-I} and \ref{maxdecaynonlinearX-II} are optimal with respect to the
decay. In particular, if $u$ is any nontrivial solution of (\ref{nls}), then the following estimate never occurs.
\begin{gather}
\label{decmaxresprincL2-1}
\forall t>0,\; \|u(t)\|_{L^r} \le C|t|^{-N\left(\frac{1}{2}-\frac{1}{r}\right)}\left(\ln|t|\right)^{-\delta},
\end{gather}
for some $r>2$ and $\delta>0.$ This is very surprising since some above results are established for solutions of some heat or
Ginzburg-Landau equations (see the beginning of Section \ref{introduction}). For example, such estimates are obtained for the
solutions $u$ of the Schrödinger equation
$$
u_t-u_{xx}+i|u|^2u=0, \quad (t,x)\in[0,\infty)\times\R,
$$
if $\|u(0)\|_X$ is sufficiently small (Theorem 1.1 of Hayashi, Kaikina and Naumkin \cite{MR1800091}). Furthermore, if
$\alpha>\frac{4}{N}$ then Theorems \ref{maxdecaynonlinearX-I} and \ref{maxdecaynonlinearX-II} are optimal with respect to the
assumption on the initial data $\vphi,$ that is $\vphi\in H^1(\R^N),$ in the sense that $H^1(\R^N)$ is the smallest functional
space in which we must take $\vphi$ to have a solution. On the other hand, when $\alpha\le\frac{4}{N},$ we have to make the
additional assumption on initial data $\vphi,$ that is $\vphi\in X.$ This request is very reasonable since this is in this
functional space that we obtain solutions of (\ref{nls}) which have a linear decay (see the references cited in Section
\ref{introduction}).
\end{rmk}

\begin{rmk}
\label{rmkmaxdecaynonlinearX-t<0}
Note that all the results of this section and Section \ref{decmaxresprincL2} hold for $t<0$ as soon as the solution $u$ is negatively global in time. Indeed, it is sufficient to apply the case $t>0$ to the solution positively global in time $\tilde u$ of (\ref{nls}) with initial data $\ovl\vphi.$ Since $\tilde u(t)=\ovl{u(-t)},$ the result for $t<0$ follows.
\end{rmk}

\section{Main results in the Lebesgue space}
\label{decmaxresprincL2}

As show the results of Section \ref{decmaxresprincX}, if we suppose a suitable asymptotic behavior of the initial value $(u(0)\in
X$ if $\alpha\le\frac{4}{N},$ $u(0)\in H^1(\R^N)$ if $\alpha>\frac{4}{N}),$ then we have a sharp lower bound. In particular, under
the hypotheses of Section \ref{decmaxresprincX}, such results do not allow estimates of type (\ref{decmaxresprincL2-1}),
for any nontrivial solution of (\ref{nls}), for some $r>2$ and $\delta>0$ (see Remark \ref{rmkmaxdecaynonlinearX-III}). In this
section, we establish some lower bounds which eventually allow estimates on the above type, only if $\alpha$ is small enough
(see Theorem \ref{maxdecaynonlinearH1-III} below). The loss of sharp estimate is compensated by a weaker assumption on $u(0),$ that
is $u(0)\in L^2(\R^N)$ if $\alpha\le\frac{4}{N}.$ As we can see, this hypothesis is optimal with respect to the integrability of
the initial data, in the sense that we make the minimal assumption on $u(0)$ to have existence of a solution. But when
$\alpha>\frac{4}{N},$ Theorems \ref{maxdecaynonlinearX-I} and \ref{maxdecaynonlinearX-II} are optimal with respect to the lower
bound and to the assumption on $u(0).$ So we only have to consider the case $\alpha\le\frac{4}{N}.$ On the other hand, if
$\alpha>\frac{4}{N+2}$ $(\alpha>2$ if $N=1),$ then the sharp estimate still holds (see Theorems \ref{maxdecaynonlinearH1-I} and
\ref{maxdecaynonlinearH1-II} below). However, we have to make an additional decay assumption on the solution $u$ $(u$ must satisfy
(\ref{maxdecaynonlinearH1-II1})).

\begin{thm}
\label{maxdecaynonlinearH1-I}
Let $\lambda\in\R\setminus\{0\},$ $\dfrac{4}{N+2}<\alpha\le\dfrac{4}{N}$ $(2<\alpha\le4$ if $N=1),$ $\vphi\in L^2(\R^N)$ and $u$
be the corresponding solution of $(\ref{nls}).$ If $\alpha=\dfrac{4}{N}$ then assume further that $u$ is positively global in time.
Suppose that for every $r\in\left[2,\frac{2N}{N-2}\right)$ $(r\in[2,\infty)$ if $N=1),$
\begin{gather}
\label{maxdecaynonlinearH1-II1}
a.e.\; t>0,\; \|u(t)\|_{L^r} \le C|t|^{-N\left(\frac{1}{2}-\frac{1}{r}\right)}.
\end{gather}
Then we have for every $r\in[2,\infty],$
\begin{gather*}
\vliminf_{t\to\infty}|t|^{N\left(\frac{1}{2}-\frac{1}{r}\right)}\|u(t)\|_{L^r}>0.
\end{gather*}
\end{thm}

\begin{thm}
\label{maxdecaynonlinearH1-II}
Let $\lambda\in\R\setminus\{0\},$ let $\alpha_0=\frac{-(N-2)+\sqrt{N^2+12N+4}}{2N},$ $\alpha_0<\alpha\le\dfrac{4}{N},$ $\vphi\in
L^2(\R^N)$ and
$u$ be the corresponding solution of $(\ref{nls}).$ If $\alpha=\dfrac{4}{N}$ then assume further that $u$ is positively global in
time. If there exists $r\in[\alpha+2,\infty]$ such that $u$ satisfies $(\ref{maxdecaynonlinearH1-II1})$ then for every
$r\in[2,\infty],$
\begin{gather*}
\vliminf_{t\to\infty}|t|^{N\left(\frac{1}{2}-\frac{1}{r}\right)}\|u(t)\|_{L^r}>0.
\end{gather*}
\end{thm}

\begin{rmk}
\label{rmkdeflineardecay}
When $\vphi\in L^2(\R^N),$ the condition (\ref{maxdecaynonlinearH1-II1}) makes sense. Indeed, we have by the Strichartz' estimates
that $u\in L^q_\loc([0,\infty);L^r(\R^N)),$ for every admissible pair $(q,r).$ This yields, $u(t)\in L^r(\R^N),$ for almost
every $t>0$ and for all $r\in\left[2,\frac{2N}{N-2}\right]$ $(r\in[2,\infty)$ if $N=2,$ $r\in[2,\infty]$ if $N=1).$
\end{rmk}

\begin{rmk}
\label{rmkmaxdecaynonlinearH1-I2}
As shows Lemma \ref{lemmaxdecaynonlinear1}, Theorem \ref{maxdecaynonlinearH1-II} has less restrictive assumptions than
Theorem \ref{maxdecaynonlinearH1-I} when $\alpha>\alpha_0.$ Indeed, we do not have to suppose that $u$ satisfies
$(\ref{maxdecaynonlinearH1-II1})$ for all $r.$ We may only assume that it is satisfied for $r=\alpha+2.$ Furthermore,
estimates of type (\ref{decmaxresprincL2-1}) do not occur. Finally, Theorem \ref{maxdecaynonlinearH1-II} can be extended for
$\alpha=\alpha_0$ in the following sense. If there exists $r\in[\alpha_0+2,\infty]$ and $\eps>0$ such that for almost
every $t>0,$
\begin{gather}
\label{rmkmaxdecaynonlinearH1-I2-1}
\|u(t)\|_{L^r} \le C|t|^{-N\left(\frac{1}{2}-\frac{1}{r}\right)-\eps},
\end{gather}
then for all $t\in\R,$ $u(t)\equiv0.$ See the proof of Theorem \ref{maxdecaynonlinearH1-II} for the justification.
\end{rmk}

\begin{rmk}
\label{maxdecaynonlinearH1-III}
In the case where $\alpha\le\frac{4}{N+2}$ $(\alpha\le1$ if $N=1),$ we have the following result. Let $\lambda\in\R\setminus\{0\},$
$0<\alpha\le\dfrac{4}{N+2}$ $(0<\alpha\le1$ if $N=1),$ $\vphi\in L^2(\R^N)$ and $u$ be the corresponding solution of
$(\ref{nls}).$ Let $r\in[2,\infty].$ If there exists $\eps>0$ such that
\begin{gather}
 \begin{cases}
  \label{maxdecaynonlinearH1-III1}
   \|u(t)\|_{L^\frac{2N}{N-2}} \le C|t|^{-\left(1+\frac{4-\alpha(N+2)}{4\alpha}+\eps\right)},
                                                                                   & if \;\;N\ge3, \smallskip \\
   \|u(t)\|_{L^r} \le C|t|^{-2\left(\frac{1}{2}-\frac{1}{r}\right)-\frac{r(1-\alpha)+\alpha}{r\alpha}-\eps},
                                                                                   & if \;\; N=2,  \medskip   \\
   \|u(t)\|_{L^\infty} \le C|t|^{-\left(\frac{1}{2}+\frac{2-\alpha}{2\alpha}+\eps\right)},
                                                                                   & if \;\; N=1,
 \end{cases}
\end{gather}
for almost every $t>0,$ then for all $t\in\R,$ $u(t)\equiv0.$ See the proof of Theorem \ref{maxdecaynonlinearH1-II} for the
justification. When $N\ge3$ and $\alpha=\dfrac{4}{N+2},$ the result is the same if we have (\ref{maxdecaynonlinearH1-III1}) for
some $r\in\left[\frac{2N}{N-2},\infty\right].$ Indeed, by Lemma \ref{lemmaxdecaynonlinear1} this hypothesis leads to
(\ref{maxdecaynonlinearH1-III1}) with $r=\frac{2N}{N-2}.$ When $\alpha<\frac{4}{N+2}$ $(\alpha\le1$ if $N=1),$ estimate
(\ref{maxdecaynonlinearH1-III1}) is a very strong assumption since it implies that $u$ decays faster than the solution of the
linear equation. Furthermore, there is a gap between the admissible and the non-admissible powers of decay
(compare (\ref{maxdecaynonlinearH1-II1}) with (\ref{maxdecaynonlinearH1-III1})).
\end{rmk}

\section{Estimates at infinity}
\label{decmaxestiminfini}

\begin{prop}
\label{propscattering}
Let $\lambda\in\R\setminus\{0\},$ $m\in\{0;1\},$ $0<\alpha\le\dfrac{4}{N-2m}$ $(0<\alpha<\infty$ if $N=m=1$ and
$0<\alpha<\dfrac{4}{N-2}$ if $N\ge2$ and $m=1),$ $\vphi\in H^m(\R^N)$ and $u\in C((-T_*,T^*);H^m(\R^N))$ be the unique
corresponding solution of $(\ref{nls}).$ Assume that $T^*=\infty.$ If there exist $t_0\ge0$ and $(\gamma,\rho)$ an
admissible pair with $\frac{\gamma\alpha}{\gamma-2}<\infty$ and $2\le\frac{\rho\alpha} {\rho-2}\le\frac{2N}{N-2}$
$(2\le\frac{\rho\alpha}{\rho-2}<\infty$ if $N=2,$ $2\le\frac{\rho\alpha}{\rho-2}\le\infty$ if $N=1)$ such that $u\in
L^\frac{\gamma\alpha}{\gamma-2}((t_0,\infty);L^\frac{\rho\alpha}{\rho-2}(\R^N)),$ then the following properties hold.
\begin{enumerate}
\item
\label{propscattering1+}
For every admissible pair $(q,r),$ $u\in L^q((0,\infty);W^{m,r}(\R^N)),$
\item
\label{propscattering2+}
There exists $u_+\in H^m(\R^N)$ such that $\vlim_{t\to\infty}\|T(-t)u(t)-u_+\|_{H^m}=0.$
\end{enumerate}
A similar result holds for $t<0.$
\end{prop}

\begin{proof*}
By remark \ref{rmkmaxdecaynonlinearX-t<0}, we only have to show the case $t>0.$
We proceed in 2 steps. Set $f(u)=\lambda|u|^\alpha u.$ \\
{\bf Step 1.} $f(u)\in L^{\gamma'}((t_0,\infty);W^{m,\rho'}(\R^N)).$ \\
We first show that $u\in L^\gamma((0,\infty);W^{m,\rho}(\R^N)).$ We already know that $u\in L^q_\loc([0,\infty);W^{m,r}(\R^N)),$
for every admissible pair $(q,r).$ We have the following integral equation.
$$
\forall S\ge0,\; \forall t\ge0,\;
u(t)=T(t-S)u(S)+i\dsp\vint_S^tT(t-s)f(u(s))ds.
$$
So we have by the Hölder's inequality (applied in space-time) and Strichartz' estimates,
\begin{gather}
\label{vproofpropscattering1}
\|f(u)\|_{L^{\gamma'}((t_0,t);W^{m,\rho'})}\le C \|u\|_{L^{\frac{\gamma\alpha}{\gamma-2}}
     ((t_0,\infty);L^{\frac{\rho\alpha}{\rho-2}})}^\alpha\|u\|_ {L^\gamma((0,t);W^{m,\rho})}, \\
\label{vproofpropscattering2}
\|u\|_{L^\gamma((S,t);W^{m,\rho})} \le C + C_0 \|u\|_{L^{\frac{\gamma\alpha}{\gamma-2}}((S,\infty);L^{\frac{\rho\alpha}
{\rho-2}})}^\alpha\|u\|_{L^\gamma((S,t);W^{m,\rho})},
\end{gather}
for every $t_0\le S<t<\infty.$ Since $u\in L^\frac{\gamma\alpha}{\gamma-2}((t_0,\infty);L^\frac{\rho\alpha}{\rho-2}(\R^N)),$ there
exists $S_0>t_0$ large enough such that $C_0\|u\|_{L^{\frac{\gamma\alpha}{\gamma-2}}
((S_0,\infty);L^{\frac{\rho\alpha}{\rho-2}})}^\alpha\le1/2,$ where $C_0$ is the constant in (\ref{vproofpropscattering2}). So with
(\ref{vproofpropscattering2}), we obtain $\|u\|_{L^\gamma((S_0,t);W^{m,\rho})} \le 2C,$ for every $t>S_0.$ It follows
that $\|u\|_{L^\gamma((S_0,\infty);W^{m,\rho})} \le 2C$ and so we have $u\in L^\gamma((0,\infty);W^{m,\rho}(\R^N)).$ Hence the
result by letting $t\nearrow\infty$ in (\ref{vproofpropscattering1}).  \\
{\bf Step 2.} Conclusion. \\
By Step 1 and Strichartz' estimates, $u\in L^q((0,\infty);W^{m,r}(\R^N)),$ for every admissible pair $(q,r).$ Then
\ref{propscattering1+} follows. From the Strichartz' estimates and by the fact that $T(t)$ is an isometry on $H^m(\R^N),$ we
obtain for every $\tau>t>t_0,$
$$
\|T(-t)u(t)-T(-\tau)u(\tau)\|_{H^m}\le C\|f(u)\|_{L^{\gamma'}((t,\tau);W^{m,\rho'})}\xrightarrow{t,\tau\to\infty}0,
$$
by Step 1. Hence \ref{propscattering2+}. This concludes the proof.
\end{proof*}

\begin{rmk}
Note that by assumption, one always has $\frac{\gamma\alpha}{\gamma-2}>0.$ However, it may happen that
$\frac{\gamma\alpha}{\gamma-2}<1.$ This is clearly not a problem since the above proof still holds and that we do not use the
triangular inequality.
\end{rmk}

\begin{lem}
\label{lemmaxdecaynonlinear1}
Let $\lambda\in\R,$ $m\in\{0;1\},$ $0\le\alpha\le\dfrac{4}{N-2m}$ $(0\le\alpha<\infty$ if $N=m=1$ and $0\le\alpha<\dfrac{4}{N-2}$
if $N\ge2$ and $m=1),$ $\vphi\in H^m(\R^N)$ and $u\in C((-T_*,T^*);H^m(\R^N))$ be the corresponding solution of $(\ref{nls}).$
Assume that $T^*=\infty.$ If there exist $r\in(2,\infty],$ $\eps\ge0$ and a constant $C=C(t)>0$ such that $u(t)$ satisfies
$(\ref{rmkmaxdecaynonlinearH1-I2-1})$ for some $t>0,$ then for every $\rho\in(2,r],$ there exist $\eps(\rho)\ge0$ and $C_0(t)>0$
such that
\begin{gather}
 \label{lemmaxdecaynonlinear11a}
  \|u(t)\|_{L^\rho} \le C_0(t)t^{-N\left(\frac{1}{2}-\frac{1}{\rho}\right)-\eps(\rho)},
\end{gather}
where the function $\rho\longmapsto\eps(\rho)$ is continuous from $(2,r]$ to $[0,\infty)$ and satisfies $\eps(\rho)>0 \iff \eps>0.$
If
$C$ is independent on $t$ then $C_0$ is also independent on $t.$ Finally, if $(\ref{rmkmaxdecaynonlinearH1-I2-1})$ is satisfied for
every $t>0$ then $(\ref{lemmaxdecaynonlinear11a})$ is satisfied for every $t>0,$ and if $\vliminf_{t\to\infty}C(t)=0$ then
$\vliminf_{t\to\infty}C_0(t)=0.$
\end{lem}

\begin{proof*}
Let $\rho\in(2,r].$ Set $\theta=\dfrac{r}{\rho}\dfrac{\rho-2}{r-2},$ $\eps(\rho)=\eps\theta$ and $C_0(t)=C(t)^\theta.$ Then
$\theta\in(0,1]$ and $\theta$ satisfies $\dfrac{1}{\rho}=\dfrac{1-\theta}{2}+\dfrac{\theta}{r}.$ By Hölder's inequality and
conservation of charge, we obtain
$$
\|u(t)\|_{L^\rho}\le\|u(t)\|_{L^2}^{1-\theta}\|u(t)\|_{L^r}^\theta
\le C(t)^\theta|t|^{-N(\frac{1}{2}-\frac{1}{2})(1-\theta)-N(\frac{1}{2}-\frac{1}{r})\theta-\eps\theta}
\le C_0(t)|t|^{-N(\frac{1}{2}-\frac{1}{\rho})-\eps(\rho)}.
$$
Hence the result.
\end{proof*}

\begin{lem}
\label{lemmaxdecaynonlinear2}
Let $\lambda\in\R\setminus\{0\},$ $m\in\{0;1\},$ $0<\alpha\le\dfrac{4}{N-2m}$ $(0<\alpha<\infty$ if $N=m=1$ and
$0<\alpha<\dfrac{4}{N-2}$ if $N\ge2$ and $m=1),$ $\vphi\in H^m(\R^N)$ and $u\in C((-T_*,T^*);H^m(\R^N))$ be the corresponding solution of $(\ref{nls}).$ Assume that $T^*=\infty.$ If $u$ satisfies $(\ref{maxdecaynonlinearH1-II1})$ for every $r\in\left[2,\frac{2N}{N-2}\right)$ $(r\in[2,\infty)$ if $N=1)$  and if $\alpha>\frac{4}{N+2}$ $(\alpha>2$ if $N=1),$ then there exists an admissible pair $(\gamma,\rho)$ with $1<\frac{\gamma\alpha}{\gamma-2}<\infty$ and $2<\frac{\rho\alpha} {\rho-2}<\frac{2N}{N-2}$ $(2<\frac{\rho\alpha}{\rho-2}<\infty$ if $N=1)$ such that $u\in
L^\frac{\gamma\alpha} {\gamma-2}((1,\infty);L^\frac{\rho\alpha}{\rho-2}(\R^N)).$
\end{lem}

\begin{proof*}
We distinguish 3 cases\:: $N\ge3,$ $N=2$ and $N=1.$ \smallskip \\
{\it Case N$\:\ge$3.} Set $\rho_*=\frac{4N}{2N-\alpha(N-2)}.$ Since $0<\alpha<\frac{4}{N-2}$ then $2<\rho_*<\frac{2N}{N-2}.$ Let
$\gamma_*>2$ be such that $(\gamma_*,\rho_*)$ is an admissible pair. For this choice of $\rho_*,$ we have
$\frac{\rho_*\alpha}{\rho_*-2}=\frac{2N}{N-2}$ and $\frac{\gamma_*}{\gamma_*-2} =\frac{4}{4-\alpha(N-2)}.$ When
$\alpha<\frac{N+2}{N},$ we have $\rho_*<\frac{2N}{(N+2)-N\alpha} \iff \alpha>\frac{4}{N+2}.$ Let $\rho>\rho_*,$ $\rho$
sufficiently close to $\rho_*$ to have $\frac{\rho\alpha}{\rho-2}>2.$ If $\alpha<\frac{N+2}{N},$ then we also choose
$\rho<\frac{2N}{(N+2)-N\alpha}.$ Since $\rho>\rho_*$ then $\frac{\rho\alpha}{\rho-2}<\frac{2N}{N-2}$ and so there exists
$\gamma>2$ such that $(\gamma,\rho)$ is an admissible pair. Then $\frac{\gamma}{\gamma-2}=\frac{2\rho}{2N-\rho(N-2)}.$
By (\ref{maxdecaynonlinearH1-II1}) we have,
\begin{gather*}
    \|u\|_{L^\frac{\gamma\alpha}{\gamma-2}((1,\infty);L^\frac{\rho\alpha}{\rho-2})}^\frac{\gamma\alpha}{\gamma-2}
 =  \vint_1^\infty\|u(t)\|_{L^\frac{\rho\alpha}{\rho-2}}^\frac{\gamma\alpha}{\gamma-2}dt
\le C\vint_1^\infty t^{-N\frac{\rho\alpha-2(\rho-2)}{2N-\rho(N-2)}}dt < \infty.
\end{gather*}
Indeed, if $\alpha<\frac{N+2}{N}$ then $N\frac{\rho\alpha-2(\rho-2)}{2N-\rho(N-2)}>1 \iff \rho<\frac{2N}{(N+2)-N\alpha}$
and if $\alpha\ge\frac{N+2}{N}$ then we always have $N\frac{\rho\alpha-2(\rho-2)}{2N-\rho(N-2)}>1.$ So, for this choice of
$(\gamma,\rho),$ $u\in L^\frac{\gamma\alpha}{\gamma-2}((1,\infty);L^\frac{\rho\alpha}{\rho-2}(\R^N)).$ \\
{\it Case N=2.} Since $\alpha>1$ is fixed, we can choose $\rho>2$ sufficiently close to 2 to have
$\alpha>\frac{2(\rho-1)}{\rho}.$ In particular, this implies that $\frac{\rho\alpha}{\rho-2}>2.$ Moreover, $\frac{\gamma}{\gamma-2}
=\frac{\rho}{2}$ where $\gamma>2$ is such that $(\gamma,\rho)$ is an admissible pair. By (\ref{maxdecaynonlinearH1-II1}) we have,
\begin{gather*}
    \|u\|_{L^\frac{\gamma\alpha}{\gamma-2}((1,\infty);L^\frac{\rho\alpha}{\rho-2})}^\frac{\gamma\alpha}{\gamma-2}
 =  \vint_1^\infty\|u(t)\|_{L^\frac{\rho\alpha}{\rho-2}}^\frac{\rho\alpha}{2}dt
\le C\vint_1^\infty t^{-\frac{\rho\alpha-2(\rho-2)}{2}}dt < \infty,
\end{gather*}
since $\frac{\rho\alpha-2(\rho-2)}{2}>1 \iff \alpha>\frac{2(\rho-1)}{\rho}.$ So $u\in L^\frac{\gamma\alpha}{\gamma-2}
((1,\infty);L^\frac{\rho\alpha} {\rho-2}(\R^2))$ for this choice of $(\gamma,\rho).$ \\
{\it Case N=1.} Since $\alpha>2$ is fixed, we can choose $\rho>2$ sufficiently close to 2 to have $\alpha>\frac{3\rho-2}{\rho}.$
In particular, this implies that $\frac{\rho\alpha}{\rho-2}>2.$ Moreover, $\frac{\gamma}{\gamma-2}=\frac{2\rho}{\rho+2}$ where
$\gamma>2$ is such that $(\gamma,\rho)$ is an admissible pair. By (\ref{maxdecaynonlinearH1-II1}) we have,
\begin{gather*}
    \|u\|_{L^\frac{\gamma\alpha}{\gamma-2}((1,\infty);L^\frac{\rho\alpha}{\rho-2})}^\frac{\gamma\alpha}{\gamma-2}
 =  \vint_1^\infty\|u(t)\|_{L^\frac{\rho\alpha}{\rho-2}}^\frac{2\rho\alpha}{\rho+2} dt
\le C\vint_1^\infty t^{-\frac{\rho\alpha-2(\rho-2)}{\rho+2}}dt < \infty,
\end{gather*}
since $\frac{\rho\alpha-2(\rho-2)}{\rho+2}>1 \iff \alpha>\frac{3\rho-2}{\rho}.$ So for this choice of $(\gamma,\rho),$ $u\in
L^\frac{\gamma\alpha}{\gamma-2}((0,\infty);L^\frac{\rho\alpha}{\rho-2}(\R)).$
\medskip
\end{proof*}

As seen in Section \ref{introduction}, the crux of the proof of results of this paper is based on the following lemma.

\begin{lem}
\label{lemmaxdecaynonlinear3}
Let $\lambda\in\R\setminus\{0\},$ $m\in\{0;1\},$ $0\le\alpha\le\dfrac{4}{N-2m}$ $(0\le\alpha<\infty$ if $N=m=1$ and
$0\le\alpha<\dfrac{4}{N-2}$ if $N\ge2$ and $m=1),$ $\vphi\in H^m(\R^N)$ and $u\in C((-T_*,T^*);H^m(\R^N))$ be the corresponding solution of $(\ref{nls}).$ Assume that $T^*=\infty.$ If $\vphi\not\equiv0$ and if there exists $u_+\in L^2(\R^N)$
such that $\vlim_{t\to\infty}\|T(-t)u(t)-u_+\|_{L^2}=0,$ then for every $r\in[2,\infty],$
\begin{gather*}
\vliminf_{t\to\infty}|t|^{N\left(\frac{1}{2}-\frac{1}{r}\right)}\|u(t)\|_{L^r}>0.
\end{gather*}
\end{lem}

The proof of Lemma \ref{lemmaxdecaynonlinear3} is based on the pseudo-conformal transformation.

For every positively global solution $u$ of (\ref{nls}) with initial data $\vphi\in L^2(\R^N),$ we define the function
$v\in C([0,1);L^2(\R^N))$ by
\begin{gather}
\label{vproofdecmaxresprincX1}
\forall t\in[0,1),\; a.e.\; x\in\R^N,\;
v(t,x)=(1-t)^{-\frac{N}{2}}u\left(\frac{t}{1-t},\frac{x}{1-t}\right)e^{-i\frac{|x|^2}{4(1-t)}}.
\end{gather}
A straightforward calculation gives for every $p\in[1,\infty]$ and for all $t\in[0,1),$
\begin{align}
 \label{vproofdecmaxresprincX2}
  \|v(t)\|_{L^p} & =(1-t)^{-N\left(\frac{1}{2}-\frac{1}{p}\right)}\left\|u\left(\frac{t}{1-t}\right)\right\|_{L^p}, \\
 \label{vproofdecmaxresprincX3}
  \|v(t)\|_{L^2} & =\|\vphi\|_{L^2},
\end{align}
where the last identity comes from (\ref{vproofdecmaxresprincX2}) and from conservation of charge for $u.$ Note that
(\ref{vproofdecmaxresprincX2}) makes sense as soon as $u\left(\frac{t}{1-t}\right)\in L^p(\R^N).$ When $\vphi\in X,$ we obviously
have $v\in C([0,1);X)$ and so we may define for all $t\in[0,1),$
\begin{gather*}
E_1(t)=\frac{1}{2}(1-t)^\frac{4-N\alpha}{2}\|\nabla v(t)\|_{L^2}^2
      -\frac{\lambda}{\alpha+2}\|v(t)\|_{L^{\alpha+2}}^{\alpha+2}, \\
E_2(t)=\frac{1}{8}\|(x+2i(1-t)\nabla)v(t)\|_{L^2}^2-\frac{\lambda}
       {\alpha+2}(1-t)^\frac{N\alpha}{2}\|v(t)\|_{L^{\alpha+2}}^{\alpha+2}.
\end{gather*}
Then for all $t\in[0,1),$
\begin{gather}
\label{vproofdecmaxresprincX4}
\frac{d}{dt}E_1(t)=\frac{N\alpha-4}{4}(1-t)^\frac{2-N\alpha}{2}\|\nabla v(t)\|_{L^2}^2, \\
\label{vproofdecmaxresprincX5}
\frac{d}{dt}E_2(t)=0.
\end{gather}
For the proof, see Proposition 3.8 and formulas (3.20) and (3.21) of Cazenave and Weissler \cite{MR93d:35150}.
\medskip

\begin{vproof}{of Lemma \ref{lemmaxdecaynonlinear3}.}
We argue by contradiction. Let $v\in C([0,1);L^2(\R^N))$ be the function defined by (\ref{vproofdecmaxresprincX1}). Assume that
there exists $r\ge2$ such that
$$
\vliminf_{t\to\infty}t^{N\left(\frac{1}{2}-\frac{1}{r}\right)}\|u(t)\|_{L^r}=0.
$$
Then, we shall show that $\vphi\equiv0.$ \\
By conservation of charge and Lemma \ref{lemmaxdecaynonlinear1}, we may assume that $2<r<\frac{2N}{N-2}$ $(2<r<\infty$ if $N=1).$
Since $u(t)\in L^r(\R^N)$ for almost every $t>0,$ it follows that $v(t)\in L^r(\R^N),$ for almost every $t\in(0,1).$ By
(\ref{vproofdecmaxresprincX2}), we have
\begin{gather}
\label{vproofmaxdecaynonlinearX-I1}
\liminf_{t\nearrow1}\|v(t)\|_{L^r}=0.
\end{gather}
By hypothesis, $\vlim_{t\to\infty}\|T(-t)u(t)-u_+\|_{L^2}=0$ for some $u_+\in L^2(\R^N).$ From Proposition 3.14 of
Cazenave and Weissler \cite{MR93d:35150}, this implies that there
exists $w\in L^2(\R^N)$ such that
\begin{gather}
\label{vproofmaxdecaynonlinearX-I2}
\vlim_{t\nearrow1}\|v(t)-w\|_{L^2}=0.
\end{gather}
(Although Proposition 3.14 is given with $\alpha>0,$ the result still holds for $\alpha=0$ since the proof applies without any modification.) From (\ref{vproofmaxdecaynonlinearX-I1}) and (\ref{vproofmaxdecaynonlinearX-I2}) we deduce that $\vlim_{t\nearrow1}\|v(t)\|_{L^2}=0,$ from which it follows with the conservation of charge (\ref{vproofdecmaxresprincX3}), $\|\vphi\|_{L^2}=0.$ This is absurd since $\vphi\not\equiv0.$
\end{vproof}

\section{Proof of the results of Section \ref{decmaxresprincX}}
\label{vproofdecmaxresprincX}

Our strategy is the following. We show that if a solution $u$ of (\ref{nls}) has a decay rate too fast, then the corresponding
function $v$ given by the pseudo-conformal transformation must converge to 0 in a Lebesgue space $L^p(\R^N),$ for some
$2<p<\infty.$ But these functions also satisfy the conservation of charge. And by using the embedding $Y\inj L^{p'}(\R^N)$
or the existence of a strong limit for $v(t)$ in $L^2(\R^N)$ as $t\nearrow1,$ we deduce that $v(t)\equiv0,$ that is $u(t)\equiv0,$
for all $t\in\R.$

In order to show Theorems \ref{maxdecaynonlinearX-I} and \ref{maxdecaynonlinearX-II}, we split the proof in 2 cases, which are
$\alpha\le\frac{4}{N}$ and $\alpha>\frac{4}{N}.$

\begin{lem}
\label{lemvproofdecmaxresprincX}
Let $\lambda\in\R,$ $0\le\alpha\le\dfrac{4}{N},$ $\vphi\in X$ and let $u\in C((-T_*,T^*);X)$ be the corresponding solution of
$(\ref{nls}).$ If $\alpha=\dfrac{4}{N}$ then we suppose that $T^*=\infty.$ If $\vphi\not\equiv0$ then
$$
\vliminf_{t\to\infty}|t|^{N\left(\frac{1}{2}-\frac{1}{r}\right)}\|u(t)\|_{L^r}>0,
$$
for every $r\in[2,\infty]$ if $\lambda\le0,$ and for every $r\in[\alpha+2,\infty]$ if $\lambda>0.$
\end{lem}

\begin{proof*}
We argue by contraposition. Let $v\in C([0,1);X)$ be the function defined by (\ref{vproofdecmaxresprincX1}). Assume there exists
$r\ge2$ if $\lambda\le0,$ and $r\ge\alpha+2$ if $\lambda>0,$ such that
$$
\vliminf_{t\to\infty}t^{N\left(\frac{1}{2}-\frac{1}{r}\right)}\|u(t)\|_{L^r}=0.
$$
Then, we have to show that $\vphi\equiv0.$ \\
By conservation of charge, if $r=2$ then $\vphi\equiv0.$ So we may assume that $r>2.$ Furthermore, by Lemma
\ref{lemmaxdecaynonlinear1}, we also may assume that $r<\frac{2N}{N-2}$ $(r<\infty$ if $N=1)$ if $\lambda\le0$ or if $\alpha=0,$
and $r=\alpha+2,$ if $\lambda>0$ and if $\alpha>0.$ Since
$\vliminf_{t\to\infty}t^{N\left(\frac{1}{2}-\frac{1}{r}\right)}\|u(t)\|_{L^r}=0,$ it follows from (\ref{vproofdecmaxresprincX2})
that $\vliminf_{t\nearrow1}\|v(t)\|_{L^r}=0.$ Thus, there exists a sequence $(t_n)_{n\in\N}\subset(0,1)$ satisfying
$t_n\xrightarrow{n\to\infty}1$ such that
\begin{gather}
\label{lemvproofdecmaxresprincX1}
\lim_{n\to\infty}\|v(t_n)\|_{L^r}=0.
\end{gather}
If $\lambda\le0$ or if $\alpha=0$ then by (\ref{vproofdecmaxresprincX3}) and (\ref{vproofdecmaxresprincX4}), we have
$\vsup_{t\in[0,1)}(1-t)\|\nabla v(t)\|_{L^2}<\infty,$ which leads with (\ref{vproofdecmaxresprincX5}) and
(\ref{vproofdecmaxresprincX3}), $\vsup_{t\in[0,1)}\|v(t)\|_Y<\infty.$ If $\lambda>0$ and if $\alpha>0$ then by
(\ref{vproofdecmaxresprincX4}) and (\ref{lemvproofdecmaxresprincX1}), we have for all $n\in\N,$ $\|\nabla v(t_n)\|_{L^2}\le
C(1-t_n)^\frac{N\alpha-4}{4}.$ It follows that,
$$
(1-t_n)\|\nabla v(t_n)\|_{L^2}\le C(1-t_n)^\frac{N\alpha}{4}\xrightarrow{n\to\infty}0,
$$
and with (\ref{lemvproofdecmaxresprincX1}) and (\ref{vproofdecmaxresprincX5}), we deduce that
$
\vsup_{n\in\N}\|xv(t_n)\|_{L^2}<\infty.
$
This last estimate yields with (\ref{vproofdecmaxresprincX3}),
\begin{gather}
\label{lemvproofdecmaxresprincX2}
\sup_{n\in\N}\|v(t_n)\|_Y<\infty.
\end{gather}
It follows that for $\lambda\in\R$ and for $\alpha\in\left[0,\frac{4}{N}\right],$ we have (\ref{lemvproofdecmaxresprincX2}). From
(\ref{vproofdecmaxresprincX3}), Hölder's inequality, from the embedding $Y\inj L^{r'}(\R^N),$
from (\ref{lemvproofdecmaxresprincX2}) and (\ref{lemvproofdecmaxresprincX1}), we obtain
$$
\|\vphi\|_{L^2}=\|v(t_n)\|_{L^2}\le\|v(t_n)\|_{L^{r'}}^\frac{1}{2}\|v(t_n)\|_{L^r}^\frac{1}{2}
\le C\|v(t_n)\|_Y^\frac{1}{2}\|v(t_n)\|_{L^r}^\frac{1}{2}
\le C\|v(t_n)\|_{L^r}^\frac{1}{2}\xrightarrow{n\to\infty}0.
$$
So $\|\vphi\|_{L^2}=0$ which is $\vphi\equiv0.$ Hence the result.
\medskip
\end{proof*}

\begin{vproof}{of Theorem \ref{maxdecaynonlinearX-I}.}
If $\alpha\le\frac{4}{N}$ then the result comes from Lemma \ref{lemvproofdecmaxresprincX}. So we may assume that
$\alpha>\frac{4}{N}.$ Since $\lambda<0$ and $\alpha>\frac{4}{N},$ there exists $u_+\in H^1(\R^N)$ such that
$\vlim_{t\to\infty}\|T(-t)u(t)-u_+\|_{H^1}=0$ (Ginibre and Velo \cite{MR87i:35171}, Nakanishi \cite{MR1726753,MR1829982}). The
result comes from Lemma \ref{lemmaxdecaynonlinear3}.
\medskip
\end{vproof}

\begin{vproof}{of Theorem \ref{maxdecaynonlinearX-II}.}
We proceed in 4 steps. Let $v\in C([0,1);L^2(\R^N))$ be the function defined by (\ref{vproofdecmaxresprincX1}). \\
{\bf Step 1.} If $\alpha>\frac{4}{N}$ and if $\vlimsup_{t\to\infty}t^{N\left(\frac{1}{2}-\frac{1}{\alpha+2}\right)}
                 \|u(t)\|_{L^{\alpha+2}}\le C$ then there exists $u_+\in H^1(\R^N)$ such that
                 \begin{gather}
                  \label{vproofmaxdecaynonlinearX-II1}
                   \vlim_{t\to\infty}\|T(-t)u(t)-u_+\|_{H^1}=0.
                 \end{gather}
Let $q=\frac{4(\alpha+2)}{N\alpha}.$ Then $(q,\alpha+2)$ is an admissible pair. Since $u\in C([0,\infty);H^1(\R^N))$ and that
$H^1(\R^N)\inj L^{\alpha+2}(\R^N),$ then $u\in L_\loc^\frac{q\alpha}{q-2}([0,\infty);L^{\alpha+2}(\R^N)).$ Since
$\alpha>\frac{4}{N},$ then $\frac{N\alpha^2}{4-\alpha(N-2)}>1$ and it follows that
$$
    \|u\|^\frac{q\alpha}{q-2}_{L^\frac{q\alpha}{q-2}((1,\infty);L^{\alpha+2})}
 =  \vint_1^\infty\|u(t)\|_{L^{\alpha+2}}^\frac{q\alpha}{q-2}dt
\le C\vint_1^\infty t^{-\frac{N\alpha^2}{4-\alpha(N-2)}}dt < \infty.
$$
Therefore, $u\in L^\frac{q\alpha}{q-2}((0,\infty);L^{\alpha+2}(\R^N))$ and the result comes from 
Proposition \ref{propscattering}. \\
{\bf Step 2.} If $\vphi\not\equiv0$ and if $\alpha\le\frac{4}{N}$ then for all $r\ge\alpha+2,$
                   $\vliminf_{t\to\infty}t^{N\left(\frac{1}{2}-\frac{1}{r}\right)}\|u(t)\|_{L^r}>0.$ \\
The result comes from Lemma \ref{lemvproofdecmaxresprincX}. \\
{\bf Step 3.} If $\vphi\not\equiv0$ and if $\alpha>\frac{4}{N}$ then for all $r\ge\alpha+2,$
                 $\vlimsup_{t\to\infty}t^{N\left(\frac{1}{2}-\frac{1}{r}\right)}\|u(t)\|_{L^r}>0.$ \\
We argue by contraposition. Assume that there exists $r\ge\alpha+2$ such that
$$
\limsup_{t\to\infty}t^{N\left(\frac{1}{2}-\frac{1}{r}\right)}\|u(t)\|_{L^r}=0.
$$
Then, we have to show that $\vphi\equiv0.$ \\
By Lemma \ref{lemmaxdecaynonlinear1}, we may assume that $r=\alpha+2.$ Step 1 implies that there exists $u_+\in H^1(\R^N)$
satisfying (\ref{vproofmaxdecaynonlinearX-II1}). Then $\vphi\equiv0$ by Lemma \ref{lemmaxdecaynonlinear3}, which is the desired
result. \\
{\bf Step 4.} If $\vphi\not\equiv0$ and if there exists $\rho\ge\alpha+2$ such that
$
\vlimsup_{t\to\infty}t^{N\left(\frac{1}{2}-\frac{1}{\rho}\right)}\|u(t)\|_{L^\rho}<\infty
$
then
$$
\liminf_{t\to\infty}t^{N\left(\frac{1}{2}-\frac{1}{r}\right)}\|u(t)\|_{L^r}>0,
$$
for all $r\in[2,\infty].$ \\
If $\alpha=0$ then Step 2 gives the result and so we consider the case $\alpha>0.$ By Lemma \ref{lemmaxdecaynonlinear1}, we may
assume that $\rho=\alpha+2.$ We argue by contradiction. Suppose that there exists $r\ge2$ such that
$\vliminf_{t\to\infty}t^{N\left(\frac{1}{2}-\frac{1}{r}\right)}\|u(t)\|_{L^r}=0.$ Then $2<r<\frac{2N}{N-2}$ $(2<r<\infty$ if
$N=1).$ Indeed, this comes from conservation of charge and Lemma \ref{lemmaxdecaynonlinear1}. We obtain with
(\ref{vproofdecmaxresprincX2}),
\begin{gather}
 \label{vproofmaxdecaynonlinearX-II4}
  \sup_{t\in[0,1)}\|v(t)\|_{L^{\alpha+2}}<\infty, \\
 \label{vproofmaxdecaynonlinearX-II5}
  \liminf_{t\nearrow1}\|v(t)\|_{L^r}=0.
\end{gather}
Note that since $u\in C([0,\infty);H^1(\R^N))$ and that the embedding $H^1(\R^N)\inj L^r(\R^N)\cap L^{\alpha+2}(\R^N)$ holds, then
we have $v\in C([0,1);L^r(\R^N)\cap L^{\alpha+2}(\R^N)).$ \\
$\mathit{Case \; 1}:$ $0<\alpha\le\frac{4}{N}.$ \\
From (\ref{vproofdecmaxresprincX3}), (\ref{vproofdecmaxresprincX4}), (\ref{vproofdecmaxresprincX5}) and
(\ref{vproofmaxdecaynonlinearX-II4}), $\vsup_{t\in[0,1)}\|v(t)\|_Y<\infty.$ From (\ref{vproofdecmaxresprincX3}), from Hölder's
inequality and the embedding $Y\inj L^{r'}(\R^N),$ we have for all $t\in[0,1),$
$$
\|\vphi\|_{L^2}=\|v(t)\|_{L^2}\le\|v(t)\|_{L^{r'}}^\frac{1}{2}\|v(t)\|_{L^r}^\frac{1}{2}
\le C\|v(t)\|_Y^\frac{1}{2}\|v(t)\|_{L^r}^\frac{1}{2}
\le C\|v(t)\|_{L^r}^\frac{1}{2}.
$$
Thus $\|\vphi\|_{L^2}\le C\vliminf_{t\nearrow1}\|v(t)\|_{L^r}^\frac{1}{2}=0$ by (\ref{vproofmaxdecaynonlinearX-II5}) and so
$\|\vphi\|_{L^2}=0,$ which is absurd. \\
$\mathit{Case \; 2}:$ $\alpha>\frac{4}{N}.$ \\
By Step 1, there exists $u_+\in H^1(\R^N)$ satisfying (\ref{vproofmaxdecaynonlinearX-II1}), which gives $\vphi\equiv0$ by
Lemma \ref{lemmaxdecaynonlinear3}. This result being absurd, Step 4 is true. This concludes the proof.
\medskip
\end{vproof}

\begin{vproof}{of Corollary \ref{cormaxdecaynonlinearX-II}.}
By Cazenave and Weissler \cite{MR93d:35150}, we know that if $\|\vphi\|_X$ is sufficiently small, then $u$ is global in time and
there exists $u_+\in X$ such that $T(-t)u(t)\xrightarrow[t\to\infty]{X}u_+.$ Then, Lemma \ref{lemmaxdecaynonlinear3} gives the
result.
\medskip
\end{vproof}

\begin{vproof}{of Proposition \ref{maxdecaynonlinearX-III}.}
It is well-known that if $\|\vphi\|_{H^1}$ is sufficiently small then $u$ is global in time and $u\in L^{\alpha+2}
(\R;L^{\alpha+2}(\R^N))$ (Remark 7.7.6 of Cazenave \cite{caz1}). Then $T(-t)u(t)\xrightarrow[t\to\infty]{H^1(\R^N)}u_+,$ for some
$u_+\in H^1(\R^N)$ (Proposition \ref{propscattering}), and the result comes from Lemma \ref{lemmaxdecaynonlinear3}.
\end{vproof}

\section{Proof of the results of Section \ref{decmaxresprincL2}}
\label{vproofdecmaxresprincL2}

Our strategy is the same as for Section \ref{vproofdecmaxresprincX}. However, we could give an other proof as follows, without
requiring the pseudo-conformal transformation. We would show that if a solution $u$ of (\ref{nls}) had a decay rate too fast,
then $u$ would have a scattering state $u_\infty$ whose corresponding solution of the linear problem (that is (\ref{nls}) with
$\lambda=0)$ would have a decay rate of the same order of $u.$ In particular, $\alpha>\frac{2}{N}$ otherwise $u_\infty\equiv0$
(Barab \cite{MR86a:35121}, Strauss \cite{297.35062,MR83b:47074a}). This rate being too fast, we would have $u_\infty\equiv0$ (by
(\ref{maxdecaylinear})). And from conservation of charge, we would deduce that $u(t)\equiv0,$ for all $t\in\R.$ Furthermore, in
the case $N=1,$ we would have to make the additional assumption $\vphi\in X$ when $1<\alpha\le2$ (in order to apply the result of
Barab \cite{MR86a:35121}). But this case falls into the scope of Theorems \ref{maxdecaynonlinearX-I} and
\ref{maxdecaynonlinearX-II} where there is a better result. It follows that in this case, the result would not be interesting.
\medskip

\begin{vproof}{of Theorems \ref{maxdecaynonlinearH1-I} and \ref{maxdecaynonlinearH1-II} and Remarks
\ref{rmkmaxdecaynonlinearH1-I2} and \ref{maxdecaynonlinearH1-III}.}
We proceed in 2 steps. \\
{\bf Step 1.} There exists $u_+\in L^2(\R^N)$ such that $\vlim_{t\to\infty}\|T(-t)u(t)-u_+\|_{L^2}=0.$ \\
{\it Case of Theorems \ref{maxdecaynonlinearH1-I}.}
Since $u$ satisfies (\ref{maxdecaynonlinearH1-II1}) for every $r\in\left[2,\frac{2N}{N-2}\right)$ $(r\in[2,\infty)$ if $N=1),$ it
follows from Lemma \ref{lemmaxdecaynonlinear2} that there exists an admissible pair
$(\gamma,\rho)$ such that $u\in L^\frac{\gamma\alpha}{\gamma-2}((1,\infty);L^\frac{\rho\alpha}{\rho-2}(\R^N)).$ The result follows
from Proposition \ref{propscattering}. \\
{\it Case of Theorems \ref{maxdecaynonlinearH1-II} and Remark \ref{rmkmaxdecaynonlinearH1-I2}.}
Set $q=\frac{4(\alpha+2)}{N\alpha}.$ Thus $(q,\alpha+2)$ is an admissible pair. By Lemma \ref{lemmaxdecaynonlinear1}, we may
assume that $r=\alpha+2$ in (\ref{maxdecaynonlinearH1-II1}) and in (\ref{rmkmaxdecaynonlinearH1-I2-1}). Let $\eps>0$ as in
(\ref{rmkmaxdecaynonlinearH1-I2-1}) $(\eps=0$ in (\ref{maxdecaynonlinearH1-II1})). We set $\eps_0=\frac{q\alpha}{q-2}\eps.$ And
since $\frac{N\alpha^2}{4-\alpha(N-2)}\ge1 \iff \alpha\ge\alpha_0,$ it follows from (\ref{maxdecaynonlinearH1-II1}) or
(\ref{rmkmaxdecaynonlinearH1-I2-1}) that,
$$
    \|u\|^\frac{q\alpha}{q-2}_{L^\frac{q\alpha}{q-2}((1,\infty);L^{\alpha+2})}
 =  \vint_1^\infty\|u(t)\|_{L^{\alpha+2}}^\frac{q\alpha}{q-2}dt
\le C\vint_1^\infty t^{-\frac{N\alpha^2}{4-\alpha(N-2)}-\eps_0}dt < \infty.
$$
Then $u\in L^\frac{q\alpha}{q-2}((1,\infty);L^{\alpha+2}(\R^N))$ and the result comes from Proposition \ref{propscattering}. \\
{\it Case of Remark \ref{maxdecaynonlinearH1-III}.}
Let $r\ge2$ and $\eps>0$ be as in (\ref{maxdecaynonlinearH1-III1}). By conservation of charge, $r>2.$
Furthermore when $N=2,$ we may assume that $r<\infty$ (Lemma \ref{lemmaxdecaynonlinear1}). Let
$(\gamma,\rho)=\left(\frac{8}{\alpha(N-2)},\frac{4N}{2N-\alpha(N-2)}\right)$ if $N\ge3,$
$(\gamma,\rho)=\left(\frac{2r}{\alpha},\frac{2r}{r-\alpha}\right)$ if $N=2$ and $(\gamma,\rho)=(\infty,2)$ if $N=1.$ Then,
$\frac{\rho\alpha}{\rho-2}=\frac{2N}{N-2}$ if $N\ge3,$ $\frac{\rho\alpha}{\rho-2}=r$ if $N=2$ and
$\frac{\rho\alpha}{\rho-2}=\infty$ if $N=1.$ Applying (\ref{maxdecaynonlinearH1-III1}), it follows that for these choices of
$(\gamma,\rho),$ $u\in L^\frac{\gamma\alpha}{\gamma-2}((1,\infty);L^\frac{\rho\alpha}{\rho-2}(\R^N)).$ The result comes from
Proposition \ref{propscattering}. \\
{\bf Step 2.} Conclusion. \\
The result comes from Step 1 and Lemma \ref{lemmaxdecaynonlinear3}. This achieves the proof.
\bigskip
\end{vproof}

\noindent
{\large\bf Acknowledgments} \\
The author wishes to thank his thesis adviser, Professor Thierry Cazenave, for having suggested this work and for helpful advices.
\bigskip \\
{\large\bf Note added in proof.} Recently, a generalization of Theorem~\ref{maxdecaynonlinearX-I} has been
established for a large class of nonlinearities, as soon as the solution is bounded in time in $H^1_0 (\Omega).$ Unfortunately, these results do not apply in the case of $L^2-$solutions (which is the case in Section~\ref{decmaxresprincL2} of this paper). For more details, see~\cite{MR2100034}.

\baselineskip .0cm

\bibliographystyle{abbrv}
\bibliography{BiblioPaper3}

\end{document}